\documentclass[12pt]{amsart}
\usepackage{amsmath}
\usepackage{amsfonts}\usepackage{verbatim}
\usepackage{amssymb}\usepackage{marvosym}\usepackage{wasysym}
\usepackage{dictsym}
\usepackage{amstext}
\usepackage{amsbsy}
\usepackage{amsopn}\usepackage{MnSymbol}
\usepackage{amsthm}\usepackage{color}

\theoremstyle{definition}

\theoremstyle{remark}

\def\proclaim#1{\vskip0.5em\noindent{\bf #1}\it}
\def\endproclaim{\vskip0.5em\par\noindent\rm}
\def\proclaim#1{\vskip0.5em\noindent{\bf #1}\it}
\def\endproclaim{\vskip0.5em\par\noindent\rm}
\def\demo#1{\vskip0.5em\noindent{\bf #1\ }}

\def\text#1{\mbox{#1}}
\def\flushpar{\par\noindent}

\newcommand{\mapright}[1]{%
    \smash{\mathop{%
        \hbox to 1cm{\rightarrowfill}
        }
    \limits^{#1}
    }
}
\newcommand{\mapleft}[1]{%
    \smash{\mathop{%
        \hbox to 1cm{\rightarrowfill}
        }
    \limits_{#1}
    }
}

\def\e{\epsilon}
\def\a{\alpha}

\def\D{\Delta}
\def\s{\sigma}

\def\th{\theta}
\def\l{\lambda}
\def\x{\times}

\def\f{\flushpar}
\def\u{\underline}
\def\vf{\varphi}

\def\om{\omega}
\def\Om{\Omega}
\def\B{\mathcal B}

\def\({\biggl(}
\def\){\biggr)}

\def\<{\langle}
\def\>{\rangle}
\def\bul{\smallskip\f$\bullet\ \ \ $}
\def\Smi{\smallskip\f{ \Smiley\ \ \ }}

\def\pf{\smallskip\f{\it Proof}\ \ \ \ }\def\sbul{\f$\bullet\
\ \ $}\def\Smi{\smallskip\f{\Large\bf\Smiley\ \ \ }}
\def\lra{\longrightarrow}\def\Lra{\Longrightarrow}
\def\st{\text{such that}}\def\incss{\Bbb N^\Bbb N(\uparrow)}\def\xyr{\xrightarrow}

\begin{document}\title{ Conditions for rational weak mixing  }
\author{Jon. Aaronson}
\address{\ \ School of Math. Sciences, Tel Aviv University,
69978 Tel Aviv, Israel.}
\email{aaro@tau.ac.il}


\thanks{ The author's research was partially
supported by Israel Science Foundation grant No. 1114/08. The author would like to thank   the University of Surrey, for hospitality provided when this paper was written. }

\begin{abstract}We exhibit  rationally ergodic,  spectrally weakly mixing  measure preserving transformations which are not subsequence  
rationally weakly mixing and give a condition for smoothness of renewal sequences. \end{abstract}
\maketitle\markboth{Jon. Aaronson \copyright  2012}{weak mixing}
\

\subsection*{Weak mixing in infinite ergodic theory}
\

  The classical ``weak mixing theorem'' breaks down for conservative, ergodic transformations preserving infinite measures (see \S4 of \cite{ALW}).
However  (Thm. 4.7 in \cite{ALW}) for   a conservative, ergodic, measure preserving transformation $(X,\B,m,T)$ of a $\s$-finite measure space, 
the following conditions are equivalent:
\begin{align*}\tag{i}f\in L^\infty,\ \l\in\mathbb S^1,\ f\circ T=\l f\ \text{\tt a.e.}\ \Rightarrow\ f\ \text{\tt is constant a.e.}\end{align*}
\begin{align*}\tag{ii} T\x S\ \text{\tt\small is ergodic}\ \forall\ \text{\tt\small  ergodic, probability preserving}\ S;\end{align*}
\begin{align*}\tag{iii}\tfrac1n\sum_{k=0}^{n-1}|\int_Xuf\circ T^kdm|\underset{n\to\infty}\lra 0\ \forall\ u\in L^1_0,\ f\in L^\infty.\end{align*}
We'll call a  measure preserving transformation  satisfying (any one of) them {\it spectrally weakly mixing}. This is in the interest of disambiguation as other notions of 
``weak mixing'' are also considered here.
In \cite{ALW} and elsewhere ``spectral weak mixing" is called  "weak mixing".

A stronger notion was introduced in in \cite{RatWM}. The conservative, ergodic, measure preserving transformation $(X,\B,m,T)$ is called 
\bul {\it rationally weakly mixing along $\mathfrak K\subset\Bbb N$}  if
$\exists\ F\in\mathcal F_+$ so that
\begin{align*}\tag{$\bigstar$}\frac1{a_n(F)}\sum_{k=0}^{n-1}|m(A\cap T^{-k}B)-m(A)m(B)u_k(F)|&\underset{n\to\infty,\ n\in\mathfrak K}\lra\ 0\ \\ &  \forall\ A,\ B\in \B\cap F;\end{align*}
where $u_k(F):=\tfrac{m(F\cap T^{-k}F)}{m(F)^2}$;
\bul {\it rationally weakly mixing}  if it is  rationally weakly mixing along $\Bbb N$ and 
\bul {\it subsequence rationally weakly mixing}  if it is rationally weakly mixing along some $\mathfrak K\subset\Bbb N$.
\

\subsection*{Weak rational ergodicity}\ \
\

As in \cite{RatErg}, 
the conservative, ergodic, measure preserving transformation $(X,\B,m,T)$ is called {\it weakly rationally ergodic}  if $\exists$   $F\in\mathcal F_+$ so  that
\begin{align*}\tag{$\largestar$}\label{largestarK}\frac1{a_n(F)}\sum_{k=0}^{n-1}m(B\cap T^{-k}C)\underset{n\to\infty,\ n\in\mathfrak K}\lra\ m(B)m(C)\ \forall\ B,C\in\B\cap F\end{align*}
where $a_n(F):=\frac1{m(F)^2}\sum_{k=0}^{n-1}m(F\cap T^{-k}F)$.

In case $T$ is weakly rationally ergodic
\sbul the collection of sets $R(T)$ satisfying ($\largestar$)  is  a hereditary ring;
\sbul  $\exists\ a_n(T)$ (the {\it return sequence}) $\st$
$$\frac{a_n(A)}{a_n(T)}\underset{n\to\infty}\lra 1\ \forall\ A\in R(T);$$
\sbul as shown in \cite{RatWM}, if $T$ is rationally weakly mixing, then $T$ is weakly rationally ergodic and 
$$\{F\in\mathcal F_+:\ \text{($\bigstar$) holds for $F$}\}=R(T).$$

The stronger properties of rational ergodicity and bounded rational ergodicity are considered in \cite{RatErg} and \cite{BRE} respectively.

\

We show by  example in \S1 that
weak rational ergodicity and  spectral weak mixing does not imply subsequence rational weak mixing. 
The main examples are dyadic towers with super-growth sequences (defined below).  
We also give zero type examples.  See \cite{DGPS} for related examples.

\

As shown in \cite{HR}, a Markov shift is conservative, ergodic iff the associated stochastic matrix is irreducible and recurrent, 
and in this case is rationally ergodic (\cite{RatErg}). It is  spectrally weakly mixing iff its associated renewal sequences are aperiodic, and
(subsequence) rationally weakly mixing iff the associated renewal sequences are (subsequence) smooth (thm. 9.1 in \cite{RatWM}).
\

It is not known whether every aperiodic, recurrent renewal sequence  is subsequence smooth, or whether smoothness implies Orey's strong ratio limit property. See \S8 in \cite{RatWM}.
\

Smoothness of the renewal sequence $u$ would  follow e.g. from the property $\sum_{n\ge 1}|u_n-u_{n+1}|<\infty$. 
This  is known for positively recurrent, aperiodic,  renewal sequences and is  conjectured  for all aperiodic  renewal sequences
(see \cite{Kingman} \S1.6(iv)).
\

In \cite{RatErg}, we showed using \cite{GL} that an aperiodic renewal sequence with $\sum_{k=1}^nu_k$  $\a$-regularly varying with $0<\a<1$
is smooth. 
Here, in \S2, we give a sufficient condition establishing smoothness e.g. when $\sum_{k=1}^nu_k$ is $1$-regularly varying.

\

\

\section*{\S1 Examples}

\subsection*{1.1 Dyadic towers over the dyadic adding machine}
\

Let $\Om:=\{0,1\}^\Bbb N,$ and let $P\in\mathcal P(\Om)$ be the symmetric product measure: $P=\prod(\frac12,\tfrac12)$, and let $\tau:\Om\to\Om$ be the dyadic odometer
defined by
$$\tau(1,\dots,1,0,\om_{\ell+1},\dots)=(0,\dots,0,1,\om_{\ell+1},\dots)$$
where $\ell=\ell(\om):=\min\,\{n\ge 1:\ \om_n=0\}$.
\

An increasing sequence $q\in\incss$ is a {\it growth sequence} as in \cite{BRE} if $q_n>\sum_{1\le k<n}q_k$.
\

The {\it dyadic cocycle } $\varphi:\Om\to\Bbb N$ associated to the growth sequence $q\in\incss$ is  defined by
\begin{align*}\tag{\Football}\varphi(\om):=q_{\ell(\om)}-\sum_{k=1}^{\ell(\om)-1}q_k=\sum_{n\ge 1}q_n(\tau(\om)_k-\om_k).\end{align*}
\
The {\it dyadic tower} with the growth sequence $q$ is
the tower over $(\Om,\B(\Om),P,\tau)$ with height function $\vf$, namely
 $( X,\B( X), m, T)$ with
 \begin{align*}&X:=\{(x,n)\in\Om\x\Bbb N:\ 1\le \vf(x)\},\ m(A\x\{n\}):=P(A\cap [\vf\ge n]),\\ & T(x,n)=\begin{cases} & (x,n+1)\ \ \ \ \vf(x)\ge n+1;\\ & (\tau x,1)\ \ \ \ \vf(x)=n.\end{cases}
 \end{align*}

\subsubsection*{Rational ergodicity of $T$}

Recall from \cite{BRE} that $(X,\B,m,T)$ is (boundedly) rationally ergodic with return sequence
$$a_n(T)\asymp 2^{c(n)}\ \text{\rm where}\ c(n) =\min\,\{k\ge 1:\ q_k\ge n\}.$$

\subsubsection*{Weak mixing of $T$}
\

By \cite{AN}, $T$ is spectrally weakly mixing iff $G_2(q)=\{0\}$ where
$$G_2(q):=\{t\in\Bbb T:\ \sum_{n\ge 1}\|q_nt\|^2<\infty\}$$
where $\|x\|:=\min_{n\in\Bbb Z}|x-n|$.
\

By lemma 3 and theorem 2 in \cite{Parreau} (see also \cite{AHL}), if $q_{n+1}=a_nq_n+1$ where $a_n\in\Bbb N,\ \sum_n\tfrac1{a_n^2}=\infty$, then $G_2(q)=\{0\}$ and $T$ is  spectrally weakly mixing. See also \cite{AFS}.
\subsubsection*{Negation of subsequence rational weak mixing}
\

Let $q\in\incss$ be a   growth sequence. For $\e\in\mathcal E:=\{\eta\in\{-1,0,1\}^\Bbb N:\ \ \eta_n\to 0\}$, let
$$N_\e:=\sum_{k\ge 1}\e_kq_k.$$
It is easy to see that $\e\mapsto N_\e$ ($\mathcal E\to\Bbb Z$) is injective if
$q$ is a {\it super growth sequence} in the sense that $q_n>2\sum_{1\le k<n}q_k$.
\

For $\e\in\mathcal E$, we have that $N_\e>0$ iff
$\e_{\kappa_{\max}}=1$ where $\kappa_{\max}(\e):=\max\,\{k\ge 1:\e_k\ne 0\}$.
\

Write $\mathcal E^+:=\{\e\in\mathcal E:\ \e_{\kappa_{\max}}=1\}$ and $\|\e\|:=\sum_{n\ge 1}|\e_n|$.

\

We claim that
\Smi {\em A dyadic tower with a super growth sequence cannot be subsequence, rationally weakly mixing.}
\

Using the above, it is easy to construct super growth sequences $q\in\incss$ with $G_2(q)=\{0\}$
and hence with  spectrally weakly mixing dyadic towers.
\

Let $q\in\incss$ be a super growth sequence and let $(X,\B,m,T)$ be the corresponding dyadic tower as above.
\

The claim (\Smiley) will follow from
\begin{align*}\tag{\dsrailways}m(\Om\cap T^{-n}\Om)=\begin{cases}& \frac1{2^{\|\e\|}}\ \ \ n=N_\e,\ \e\in\mathcal E;\\ & 0\ \ \ \ \ \text{\rm else.}\end{cases}\end{align*}
\demo{Proof of (\dsrailways)\ \ \ {\rm (see \cite{HajKak})}}
\

 Let $\mathcal N\ge 1$ and $x\in\Om$, then $T^\mathcal N x\in\Om$ iff
 $\exists\ N\ge 1$ so that $\vf_N(x)=\mathcal N$.

By (\Football),
$$\mathcal N=\vf_N(x)=\sum_{k\ge 1}q_k(\tau^N(\om)_k-\om_k)=\sum_{k\ge 1}\e_kq_k=:N_\e$$ for some
$\e\in\mathcal E^+$.
\

We now show that
$m(\Om\cap T^{-N_\e}\Om)=\frac1{2^{\|\e\|}}.$
\

If $n_1<n_2<\dots<n_k,\ m_1<m_2<\dots m_\ell$ and
$$\e_{n_i}=1,\ \e_{m_j}=-1\ \ \&\ \ \e_n=0\ \ \text{\rm else},$$
then
$$\Om\cap T^{-N_\e}\Om=\{\om\in\Om:\ \om_{n_i}=0\ \forall\ i\ \&\ \om_{m_j}=1\ \ \forall\ j\}$$
and for $\om\in\Om\cap T^{-N_\e}\Om,\ T^{N_\e}\om=G\om$ where $G:\Om\to\Om$ is defined by
$$G(x)_k=\begin{cases}& 1-x_k\ \ \ \ \ k\in\{n_i\}_i\cup\{m_j\}_j;\\ &x_k\ \ \ \ \text{\rm else}.\end{cases}$$
Thus
$$m(\Om\cap T^{-N_\e}\Om)=m(\{\om\in\Om:\ \om_{n_i}=1\ \forall\ i\ \&\ \om_{m_j}=0\forall\ j\})=\frac1{2^{\|\e\|}}.\ \ \ \text{\text{\Checkedbox}(\dsrailways)}$$
\demo{Proof of (\Smiley)}
\

It suffices to show that
$$\varliminf_{n\to\infty}\frac1{2^{c(n)}}\sum_{k=1}^n|u_k-u_{k+q_1}|>0$$
where $u_n=u_n(\Om):=m(\Om\cap T^{-n}\Om)$.
\

 To see this, we restrict summation to
 $k=N_\e$ where $\e\in\mathcal E^+\  \&\ \e_1=0$; noting  $N_\e\le n$ iff $\kappa_{\max}\le c(n)$.
 \

 Here
$$u_{N_\e}=\frac1{2^{\|\e\|}},\ u_{N_\e+q_1}=\frac1{2^{\|\e\|+1}}$$ and
 $$u_{N_\e}-u_{N_\e+q_1}=\frac12\cdot u_{N_\e}.$$
 \

 Thus
 \begin{align*}\sum_{k=1}^n|u_k-u_{k+q_1}|&\ge
 \sum_{\e\in\mathcal E^+,\ \e_1=0,\ \kappa_{\max}\le c(n)}|u_{N_\e}-u_{N_\e+q_1}|\\ &=\frac12\sum_{\e\in\mathcal E^+,\ \e_1=0,\ \kappa_{\max}\le c(n)}u_{N_\e}\\ &=
 \frac12\sum_{\e\in\mathcal E^+,\ \e_1=0,\ \kappa_{\max}\le c(n)}\frac1{2^{\|\e\|}}.\end{align*}
 Now
\begin{align*}\sum_{\e\in\mathcal E^+,\ \e_1=0,\ \kappa_{\max}\le c(n)}\frac1{2^{\|\e\|}}&=
 \sum_{\emptyset\ne F\subset\Bbb N\cap [2,c(n)]}\ \ \sum_{\e\in\mathcal E^+,\ \text{\tt\tiny supp}\,\e=F}\frac1{2^{\#F}}.
  \end{align*}
  For fixed $F$,
  \begin{align*}\#\{\e\in\mathcal E,\ \text{\tt supp}\,\e=F\}=2^{\#F}.\end{align*} Thus
  \begin{align*}\sum_{k=1}^n|u_k-u_{k+q_1}|&\ge \frac12\sum_{\e\in\mathcal E^+,\ \e_1=0,\ \kappa_{\max}\le c(n)}\frac1{2^{\|\e\|}}\\ &=\frac12\#\{\emptyset\ne F\subset\Bbb N\cap [2,c(n)]\}\\ &\ge
  \frac18\cdot 2^{c(n)}.\ \ \ \text{\text{\Checkedbox}}
  \end{align*}
  \subsection*{1.2 A zero type example}
  \

  As in \cite{WWS} we say that a conservative, ergodic, measure preserving transformation $(X,\B,\mu,T)$ \bul {\it of zero type} if
  $m(A\cap T^{-n}B)\xrightarrow[n\to\infty]{}\ 0$ for some (and hence all) $A,B\in\B,\ 0<m(A),\ m(B)<\infty$, and   \bul {\it of positive type} otherwise.
  
  \

  Let $(X,\B,\mu,T)$ be a spectrally weakly mixing, dyadic tower with a super-growth sequence (as above). It follows from
  (\dsrailways) that $T$ is of positive type. To get a zero type example we multiply $T$ by a suitable Markov shift -- a {\tt renewal process}.
  \
  
  \subsection*{Renewal sequences and processes}
  \
  
  As in \cite{Chung} and \cite{Kingman}, a sequence $u=(u_0,u_1\dots)\in [0,1]^{\Bbb N_0}$  is a {\it renewal sequence} if $u_0=1\ \&$ there are numbers $f_1,f_2,\dots\ge 0$ satisfying the 
  {\it renewal equation}:
  $$u_n=\sum_{r=1}^nf_ru_{n-r}\ \ \forall\ n\ge 1.$$
  It follows from this  that $\sum_{k=1}^\infty f_k\le 1$ with equality if and only if $\sum_{n=0}^\infty u_n=\infty$, in which case
  $u$ is called {\it recurrent} and $f=(f_1,f_2,\dots)\in\mathcal P(\Bbb N)$ is called its (associated) {\it lifetime distribution}.
  \

 Given a recurrent renewal sequence $u$,
  there is an irreducible, recurrent, stochastic matrix $P:\Bbb N\x\Bbb N\to [0,1]$ so that $p_{1,1}^{(n)}=u_n$. The corresponding invertible,  stationary Markov shift 
  $$(Y,\mathcal C,\nu,S)= (\Bbb N^\Bbb Z,\B(\Bbb N^\Bbb Z),\nu,\text{\tt Shift})$$ 
(with $\nu$ the stationary Markov measure with transitions given by $P$) is called the  {\it renewal process} and is conservative, ergodic by \cite{HR}. 
\

The measure $\nu$ is infinite  iff $u$ is  null recurreent (i.e. $u_n\xrightarrow[n\to\infty]{}\ 0$).
\

In case the renewal sequence $u$ is {\it aperiodic} in the sense that
$$\gcd\,\{n\ge 1:\ u_n>0\}=1,$$ the one sided renewal process $(Y_+,\mathcal C_+,\nu,\s) $ is exact ($\bigcap_{n\ge 1}\s^{-n}\mathcal C_+\overset{\nu}=\{\emptyset,Y\}$) and
the two sided renewal process $(Y,\mathcal C,\nu,S) $ is spectrally weakly mixing.
\

As in \cite{RatErg}, a set $A\in\mathcal C,\ 0<\nu(A)<\infty$ is called a {\it recurrent event} for $S$ if for 
$0=n_0\le n_1\le n_2\le \dots\le n_k$, we have
$$\nu(\bigcap_{j=0}^kS^{-n_j}A)=\nu(A)\prod_{j=1}^k u_{n_j-n_{j-1}}\ \ \text{where}\ u_n:=\frac{\nu(A\cap T^{-n}A)}{\nu(A)}.$$
The sequence $u=u(A)=(u_0,u_1,\dots)$ is a renewal sequence. 
\

In particular, it follows from the Markov property that the set 
$A=[1]_0=\{x\in\Bbb N:\ x_0=1\}$ is a recurrent event for the renewal process $S$. 
  
Recall from \S1.5 of \cite{Kingman} that a  {\it Kaluza sequence} is a sequence $u=(u_0,u_1\dots)\in (0,1]^{\Bbb N_0}$ so that 
  $u_0=1\ \&\ \frac{u_{n+1}}{u_n}\uparrow$. Evidently, if Kaluza sequences are non-increasing.
  If in addition,  $\sum_{n=0}^\infty u_n=\infty$, then $u$ is an aperiodic, recurrent renewal sequence (thm. 1.8 in \cite{Kingman}). 
  
  \
  
We can now continue the construction of the zero type example.
  \ 
  
  By theorem 3.3 in \cite{ALW}, there is a  conservative, aperiodic,  ergodic, renewal process  $(Y,\mathcal C,\nu,S)$ with recurrent event $A$ 
  whose associated renewal sequence
  $u(A):=(u_0,u_1\dots)$ is a null recurrent Kaluza sequence and so that
  so that    $S\x T$ is conservative.
  \
  
  We claim 
\smallskip\f{\Large\bf\Yinyang\ \ \ } {\em The transformation $(Y\x X,\mathcal C\otimes\B,\nu\x\mu,S\x T)$ is zero-type, rationally ergodic and spectrally weakly mixing, but not subsequence, rationally weakly mixing.}
\demo{Proof}
\

The null recurrence of $u(A)$ ensures that $S\x T$ is of zero type.
  \

  To see that $S\x T$ is spectrally weakly mixing, let
  $V$ be an ergodic, probability preserving transformation. 
  It follows (from the spectral weak mixing of $T$) that $V\x T$ is conservative, ergodic. 
  Since $S$ is the natural extension of a an exact transformation (which is  mildly mixing), we have by theorem 6.7 in \cite{ALW} that
  $$(S\x T)\x V=S\x (T\x V)$$
  is conservative, ergodic. This shows that $S\x T$ is spectrally weakly mixing.
  \

  Next, we claim that $S\x T$ is  rationally ergodic with $A\x\Om\in R(S\x T)$.
 \pf It suffices to consider the one sided Markov shift $(Y_+,\mathcal C_+,\nu,\sigma)$ and show that $\s\x T$ is  rationally ergodic with $A\x\Om\in R(\s\x T)$.
  \

 Since $u=u(A)$ is Kaluza,  $v_k:=u_k-u_{k+1}>0$ and on $A\x\Om$, the transfer operator is given by
 \begin{align*}\widehat{\sigma\x T}^n(1_{A\x\Om})(y,\om)=\widehat{\sigma}^n(1_{A})(y)1_{\Om}(T^{-n}\om)=u_n1_{\Om}(T^{-n}\om).\end{align*}
 
Writing $a_n(\Om):=\int_\Om S_k^{(T^{\pm 1})}(1_{\Om})dm$, we have by bounded rational ergodicity of $T$ that $\exists\ M>0$ so that
$$\|S_n^{(T^{\pm 1})}(1_\Om)\|_{L^\infty(\Om)}\le M a_n(\Om)\ \forall\ n\ge 1$$
and so, 
for $(y,\om)\in A\x\Om$,
   \begin{align*}\sum_{1\le k\le n}\widehat{\sigma\x T}^k(1_{A\x\Om})(y,\om)&=\sum_{1\le k\le n}u_k1_{\Om}(T^{-k}\om)\\ &=
   \sum_{1\le k\le n}u_k(S_k^{(T^{-1})}(1_{\Om})(\om)-S_{k-1}^{(T^{-1})}(1_{\Om})(\om))\\ &=\sum_{1\le k\le n}v_kS_k^{(T^{-1})}(1_{\Om})(\om)+u_nS_n^{(T^{-1})}(1_{\Om})(\om)
   \\ &\le M\sum_{1\le k\le n}v_ka_k(\Om)+Mu_na_n(\Om)\\ &= M\sum_{1\le k\le n}u_k\mu(\Om\cap T^{-k}\Om)+2Mu_na_n(\Om)
   \\ &\le 3M\int_{A\x\Om}S_n^{(S\x T)}(A\x\Om)d\nu dm.
   \end{align*}
   Rational ergodicity of $S\x T$ and $A\x\Om\in R(S\x T)$ follow from this.
   \

   To see that $S\x T$ is not rationally weakly mixing, consider $E,\  F\in R(S\x T)$ given by
   $$E:=A\x\Om,\ F:=A\x T^{-q_1}\Om,$$ then for $n=N_\e$ with $\e_1=0$, we have
   \begin{align*}\nu\x\mu(E\cap (S\x T)^{-n}E)&-\nu\x\mu(F\cap (S\x T)^{-n}F)\\ &=
   u_n(\mu(\Om\cap T^{-n}\Om)-\mu(\Om\cap T^{-(n+q_1)}\Om))\\ &=
   \frac12u_n\mu(\Om\cap T^{-n}\Om).
   \end{align*}
   Noting that
   \begin{align*}a_n(S\x T) &\sim\sum_{\e\in\mathcal E^+,\ N_\e\le n}u_{N_\e}\frac1{2^{\|\e\|}}\\ &\asymp
   \sum_{\e\in\mathcal E^+,\ \e_1=0,\ N_\e\le n}u_{N_\e}\frac1{2^{\|\e\|}}\end{align*}
   we see that
   \begin{align*}\sum_{1\le k\le n}|\nu\x\mu(E\cap (S\x T)^{-k}E)&-\nu\x\mu(F\cap (S\x T)^{-k}F)|\\ &=
   \sum_{1\le k\le n}u_n|\mu(\Om\cap T^{-n}\Om)-\mu(\Om\cap T^{-(n+q_1)}\Om)|\\ &\ge
   \sum_{\e\in\mathcal E^+,\ \e_1=0,\ N_\e\le n}u_{N_\e}\frac1{2^{\|\e\|}}\\ &\asymp a_n(S\x T).\ \ \text{\text{\Checkedbox}\Yinyang}
   \end{align*}

   \section*{\S2 Smoothness of renewal sequences}

Suppose that  $u=u^{(f)}=(u_0,u_1,\dots)$ is an aperiodic, recurrent, renewal  sequence with lifetime distribution $f\in\mathcal P(\mathbb N)$.
\

The renewal sequence $u$ is called {\it smooth} if
$$\sum_{k=1}^n|u_k-u_{k+1}|=o(a_n)\ \text{as}\ n\to\infty\ \ \text{where}\ a_n=a_n^{(u)}:=\sum_{k=1}^nu_k.$$
It follows from  \cite{GL} that if $u=(u_0,u_1,\dots)$ is an aperiodic, recurrent, renewal  sequence and $a^{(u)}$ is $\a$-regularly varying with $\a\in (0,1)$, then $u$ is smooth (see \cite{RatWM}). The case $\a=1$ follows from the next proposition (which is related to proposition 8.3 in \cite{RatWM}).
\

For  $u=u^{(f)}$, let $c_N:=f([N,\infty))$ and $L(N):=\sum_{k=1}^Nc_k$.
\proclaim{Proposition }\ \ \ If $u=(u_0,u_1,\dots)$ is an aperiodic, recurrent, renewal  sequence with \begin{align*}\tag{\Pointinghand}\varlimsup_{N\to\infty}\frac{Nc_N}{L(N)}<\frac1{\sqrt{5}+1},\end{align*}
 then $u$ is smooth.\endproclaim
\f{\bf Remark}\ \ If  $a^{(u)}$ is $t$-regularly varying for some $t>\frac{\sqrt{5}}{\sqrt{5}+1}$ , then by the renewal equation and Karamata theory,
$L(N)\propto\frac{N}{a^{(u)}(N)}$ is $(1-t)$-regularly varying, $Nc_N\sim (1-t)L(N)$ and (\Pointinghand) holds.
\demo{Proof}
We show first that
\begin{align*}\tag{\Football}\sum_{k=1}^\infty(u_k-u_{k+1})^2<\infty.\end{align*}
Let $u=u^{(f)}$ and let $c_N:=f([N,\infty))$,  $M(N):=\sum_{1\le n\le N} nf_n$ and $V(N):=\sum_{1\le n\le N} n^2f_n$.
\

By (\Pointinghand), $\exists\ R<\frac1{\sqrt 5+1}$ so that
$Nc_N\le RL(N))$ for large $N$. It follows that
$$Nc_N\le RL(N)=R(M(N)+NC_N)$$ whence
\begin{align*}\tag{\Bicycle}Nc_N\le \frac{R}{1-R}\cdot M(N)\end{align*}
for large $N$ where $\frac{R}{1-R}<\frac1{\sqrt 5}.$
\

In particular
$$M(N)\ \asymp\ L(N)\ \ \ \ \text{as}\  N\to\infty;$$
and
$$\sum_{n=1}^\infty\frac1{n^2M(n)^2}<\infty.$$
We'll use these to prove (\Football).

\

 By  Parseval's formula, and the renewal equation,
\begin{align*} \int_{-\pi}^\pi\frac{|\th|^2d\th}{|1-f(\th)|^2}<\infty\ \ \iff\ \sum_{n=1}^\infty(u_n-u_{n+1})^2<\infty\end{align*}
where $f(\th):=\sum_{n=1}^\infty f_ne^{in\th}$. By aperiodicity, $\sup_{\e\le|\th|\le\pi}|f(\th)|<1\ \forall\ \e>0$ whence
 (using symmetry)
 \begin{align*}\int_{-\pi}^\pi\frac{|\th|^2d\th}{|1-f(\th)|^2}<\infty\ \ \iff\ \ \int_0^{\e}\frac{\th^2d\th}{|1-f(\th)|^2}<\infty\ \text{\tt\small for some}\ \ \e>0.\end{align*}
For $|\th|<\frac{\pi}2$,
 \begin{align*}1-f(\th)&=
 \sum_{n=1}^\infty f_n(2\sin^2(\frac{n\th}2)-i\sin(n\th))\\ &=
 (\sum_{1\le n\le \frac{2}{\pi|\th|}}+ \sum_{n> \frac{2}{\pi|\th|}})f_n(2\sin^2(\frac{n\th}2)-i\sin(n\th)) \end{align*}
 and
 \begin{align*}|1-f(\th)|\ge |\sum_{1\le n\le \frac{2}{\pi|\th|}}|-| \sum_{n> \frac{2}{\pi|\th|}}|.\end{align*}
 Now
 \begin{align*}|\sum_{1\le n\le \frac{2}{\pi|\th|}}|&=\sqrt{\(\sum_{1\le n\le \frac{2}{\pi|\th|}}f_n2\sin^2(\frac{n\th}2)\)^2+\(\sum_{1\le n\le \frac{2}{\pi|\th|}}f_n\sin(n\th)\)^2} \\ &\ge
 |\sum_{1\le n\le \frac{2}{\pi|\th|}}f_n\sin(n\th)|
  \\ &\ge\frac{2}{\pi}\cdot|\th|
 \sum_{1\le n\le \frac{2}{\pi|\th|}}nf_n\\ &=
 \frac{2}{\pi}\cdot|\th| M(\frac{2}{\pi|\th|});\end{align*}
 and
 \begin{align*}| \sum_{n> \frac{2}{\pi|\th|}}|&=|\sum_{n>\frac{2}{\pi|\th|}}f_n(2\sin^2(\frac{n\th}2)-i\sin(n\th))|\\ &\le \sqrt{\(\sum_{n>\frac{2}{\pi|\th|}}f_n2\sin^2(\frac{n\th}2)\)^2+\(\sum_{n> \frac{2}{\pi|\th|}}f_n\sin(n\th)\)^2} \\ &\le \sqrt{5}c_{\frac{2}{\pi|\th|}}.\end{align*}
 By assumption $\exists\ \D>0$ so that for  $|\th|<\D$,
 $$c_{\frac{2}{\pi|\th|}}\le R'\frac{2|\th|}{\pi}M(\frac{2}{\pi|\th|})$$ where
 $R'=\frac{R}{1-R}=\frac1{\sqrt 5}(1-\eta)$ for some $\eta>0$;
 whence
 \begin{align*}|1-f(\th)| & \ge \frac{2}{\pi}\cdot|\th|M(\frac{2}{\pi|\th|})-\sqrt{5}c_{\frac{2}{\pi|\th|}}\\ &\ge
 \frac{2}{\pi|\th|}M(\frac{2}{\pi|\th|})(1-\sqrt 5 R')\\ &=\frac{2\eta}{\pi|\th|}M(\frac{2}{\pi|\th|}).\end{align*}

Let $N\ge 1,\ \D>\frac{\pi}N$, then
\begin{align*}\int_0^{\frac{\pi}N}\frac{\th^2d\th}{|1-f(\th)|^2} &\ll \int_0^{\frac{\pi}N}\frac{\th^2d\th}{(\th M(\frac{2}{\pi|\th|}))^2}
\\ &=\int_0^{\frac{\pi}N}\frac{d\th}{M(\frac{2}{\pi|\th|})^2}\\ &=
\sum_{n=N}^\infty\int_{\frac{\pi}{n+1}}^{\frac{\pi}n}\frac{d\th}{M(\frac{2}{\pi|\th|})^2}
\\ &\asymp \sum_{n=N}^\infty\frac1{n^2M(n)^2}<\infty.\ \ \ \text{\Checkedbox(\Football)}\end{align*}
To see smoothness, by assumption
$$\log L(N)\sim \sum_{k=1}^N\frac{c_k}{L(k)}\le \sum_{k=1}^N\frac{1}{\sqrt 5 k}+ O(1)=\frac{1}{\sqrt 5 }\log N+O(1)$$
whence
$$L(N)=O(N^{\frac{1}{\sqrt 5 }}).$$
Using the renewal equation, $a_u(n)\asymp\frac{n}{L(n)}$, whence
$a_u(n)\gg N^{1-\frac{1}{\sqrt 5 }}$ and
$$\frac{\sqrt{n}}{a_n}\ll \frac1{N^{\frac12-\frac{1}{\sqrt 5 }}}\xyr[n\to\infty]{}0.$$
By the Cauchy-Schwartz inequality,
$$\frac1{a_u(n)}\sum_{k=1}^n|u_k-u_{k+1}|\le \frac{\sqrt{n}}{a_u(n)}\cdot\sqrt{\sum_{k=1}^\infty(u_k-u_{k+1})^2}\xyr[n\to\infty]{}0.\  \ \text{\Checkedbox}$$
\

\end{document}